\documentclass[12pt]{article}
\usepackage{latexsym}
\usepackage{amsmath}
\usepackage{amssymb}
\usepackage{bbm}
\usepackage{color}
\usepackage[tmargin=1.25in,bmargin=1.25in,lmargin=1.1in,rmargin=1.1in]{geometry}

\usepackage{graphicx}
\usepackage{epsfig}
\usepackage{hyperref} 

\definecolor{Red}{rgb}{1,0,0}
\definecolor{halfgray}{gray}{0.55}
\definecolor{webgreen}{rgb}{0,.5,0}
\definecolor{webbrown}{rgb}{.6,0,0}
\definecolor{Maroon}{cmyk}{0, 0.87, 0.68, 0.32}
\definecolor{RoyalBlue}{cmyk}{1, 0.50, 0, 0}
\definecolor{Black}{cmyk}{0, 0, 0, 0}

\hypersetup{
    colorlinks=true, linktocpage=true, pdfstartpage=1, pdfstartview=FitV,
    breaklinks=true, pdfpagemode=UseNone, pageanchor=true, pdfpagemode=UseOutlines,
    plainpages=false, bookmarksnumbered, bookmarksopen=true, bookmarksopenlevel=1,
    hypertexnames=true, pdfhighlight=/O,
    urlcolor=webbrown, linkcolor=RoyalBlue, citecolor=webgreen,
}

\newtheorem{thm}{Theorem}[section]

\newtheorem{lem}[thm]{Lemma}

\newtheorem{remark}[thm]{Remark}

\makeatletter\@addtoreset{equation}{section}\makeatother

\newcommand {\E}{\mathbb{E}}
\renewcommand {\P}{\mathbb{P}}

\newcommand{\tr}{\operatorname{tr}}
\newcommand{\sn}{{(n)}}

\def\QED{\hfill $ \Box $\\}

\begin{document}

\title{A functional CLT for partial traces of random matrices}
\author{Jan Nagel} 
\maketitle

\begin{abstract} 
In this paper we show a functional central limit theorem for the sum of the first $\lfloor t n \rfloor$ diagonal elements of $f(Z)$ as a function in $t$, for $Z$ a random real symmetric or complex Hermitian $n\times n$ matrix. The result holds for orthogonal or unitarily invariant distributions of $Z$, in the cases when the linear eigenvalue statistic $\tr f(Z)$ satisfies a CLT. The limit process interpolates between the fluctuations of individual matrix elements as $f(Z)_{1,1}$ and of the linear eigenvalue statistic. 
It can also be seen as a functional CLT for processes of randomly weighted measures. 
\end{abstract}

%
%

\section{Introduction}

It is the purpose of this paper to add a new perspective to the central limit theorem for linear eigenvalue statistics. The main objects are the eigenvalues $\lambda_1,\dots ,\lambda_n$ of a random real symmetric or complex Hermitian matrix $Z$. 
Given a test function $f$, the linear statistic of these eigenvalues, denoted by $X_1^\sn(f)$, is $\tr(f(Z))= f(\lambda_1)+\dots +f(\lambda_n)$. For many distributions of eigenvalues and smooth enough functions we have, after centering, the convergence in distribution to a normal random variable,
\begin{align} \label{CLTtrace}
\tr(f(Z)) - \mathbb{E}[\tr(f(Z)) ] 
= \sum_{k=1}^n f(\lambda_k) - \mathbb{E}[f(\lambda_k)] \xrightarrow[n \rightarrow \infty ]{d} \mathcal{N}(0,\sigma_1^2(f))  .
\end{align} 
Over the last two decades, CLTs for linear eigenvalue statistics have grown into a hugely popular field of study within random matrix theory. To give a partial overview, the convergence in \eqref{CLTtrace} was proven for invariant matrix models or orthogonal polynomial ensembles \cite{jonsson1982some,johansson1998fluctuations,pastur2006limiting,shcherbina2008invariant,kriecherbauer2010fluctuations,
dumitriu2012global,duits2015global,breuer2017central,bekerman2018clt}, 
for general Wigner or Wishart matrices  
\cite{bai2008clt,lytova2009central,shcherbina2011wigner,najim2016gaussian,bai2019central}, 
for matrices of compact groups 
\cite{johansson1997random,soshnikov2000central}, 
and for non-Hermitian matrices 
\cite{rider2006gaussian,nourdin2010universal}. 
Comparing \eqref{CLTtrace} with classical CLTs, for example for sums of independent random variables, it is highly remarkable that there is no additional scaling factor $n^{-1/2}$. This phenomenon is usually attributed to the strong dependence structure of the eigenvalues. Indeed, the classical orthogonal polynomial ensembles have a joint eigenvalue density involving the Vandermonde determinant $\Delta(\lambda) = \prod_{i<j} |\lambda_i-\lambda_j|$, which leads to a repulsion of eigenvalues. It was shown however by \cite{costin1995gaussian,soshnikov2002gaussian} that in general the variance of the linear eigenvalue statistic does not remain bounded for non-smooth test functions $f$.

One sees a very different picture when, instead of the trace, we consider the fluctuations of an individual matrix element $f(Z)_{1,1}$. Limit theorems for such entries have been considered by 
\cite{lytova2009fluctuations,lytova2011non,pizzo2012fluctuations,o2013fluctuations}. 
The random variable $f(Z)_{1,1}$ depends not only on the distribution of the eigenvalues, but also on the eigenvectors. We will assume the matrix of eigenvectors to be Haar distributed on the orthogonal group (for real $Z$) or on the unitary group (for complex $Z$), and to be independent of the eigenvalues. This is satisfied for the prominent case of unitarily invariant ensembles (see Section \ref{sec:results1}). Then the central limit theorem takes the form
\begin{align}\label{CLTentry}
\sqrt{n}(f(Z)_{1,1} - \mathbb{E} [ f(Z)_{1,1}] ) \xrightarrow[n \rightarrow \infty ]{d}  \mathcal{N}(0,\sigma_0^2(f)) . 
\end{align}
Unlike for the full trace, an additional scaling is necessary. Although one might expect $f(Z)_{1,1}$ to scale as $\frac{1}{n} \tr f(Z)$, the fluctuations of the former random variable are much larger.  
We remark that in our setting the convergence \eqref{CLTentry} is in fact a consequence of \eqref{CLTtrace} (see Theorem \ref{thm:spectralCLT}).

In this paper, we show that we can in some sense interpolate between the different CLTs in \eqref{CLTentry} and \eqref{CLTtrace} by summing a varying number of diagonal elements. The main object of interest is thus the partial trace $X_t^\sn(f)$, defined by 
\begin{align}\label{defpartialtrace}
X_t^\sn(f) = \sum_{i=1}^{\lfloor tn\rfloor } f(Z)_{i,i} = \sum_{k=1}^n w_{k,t}^\sn f(\lambda_k) , 
\end{align}
which is a weighted version of the linear eigenvalue statistic $\tr f(Z)$, where the weights $w_{k,t}^\sn$ are norms of projections of the eigenvectors (see \eqref{defweights}). 
In our main result, Theorem \ref{thm:main}, we show that in a setting where the convergence \eqref{CLTtrace} of the linear eigenvalue statistic holds, the process 
\begin{align} \label{defpartialprocess}
\big( X_t^\sn(f)- \mathbb{E}[X_t^\sn(f)] \big)_{t\in [0,1]}
\end{align}
converges as $n\to \infty$ in distribution to a centered Gaussian process. The variance of the limit process at time $t$ is given by 
\begin{align}
(t -t^2) \sigma_0^2 (f)  + t^2 \sigma_1^2(f) = t \big[ (1-t) \sigma_0^2(f) + t \sigma_1^2(f)\big] . 
\end{align}
That is, the fluctuations interpolate between the limit variance of the CLTs in \eqref{CLTentry} and \eqref{CLTtrace} and, unless $\sigma_0^2 (f)=\sigma_1^2 (f)$, the limit is not a Brownian motion.

A core argument in the proof is the independence of eigenvalues and eigenvectors. Assuming a convergence as in \eqref{CLTtrace}, the main task is then to handle the fluctuations induced by summing a varying number of entries of the eigenvector matrix. 
The main ingredient for this is a functional limit theorem for sums over subblocks of Haar distributed matrices proven by \cite{donati2012truncations,beffara2014bridges}. This result itself relies on a powerful theorem of \cite{mingo2007second}, allowing to evaluate higher order cumulants for entries of Haar matrices. Our strategy also allows us to prove a functional CLT for \eqref{defpartialprocess}, when instead of the mean $\mathbb{E}[X_t^\sn(f)]$, one centers by the expectation conditioned on the eigenvalues. The result is the \emph{quenched} convergence in Theorem \ref{thm:mainbridge}, which gives a convergence in distribution under the law of the eigenvector matrices, valid for almost all (sequences of) eigenvalues. With this centering, the limit process is a Brownian bridge. This also shows that the results are not restricted to the random matrix setting, but could also be viewed in the framework of randomly weighted sums, when the weight are coming from Haar distributed matrices as in \eqref{defpartialtrace}. For example, the functional CLT of Theorem \ref{thm:mainbridge} is also true for deterministic sequences $\lambda_i$ or more general point processes, see Remark \ref{rem:quenched}.  

Convergence of partial traces has been considered before in a couple of papers for particular distributions of random matrices. If $Z$ is unitary and $f$ the identity, a functional limit theorem for the partial trace has been proven in \cite{d2000invariance}. A more general way of summing entries of unitary matrices was considered in \cite{d2003brownian}. In \cite{rains1998normal}, real symmetric matrices are considered and the statement of Theorem \ref{thm:mainbridge} is proven under a strong moment condition on the $\lambda_i$, using zonal polynomials. Using the arguments of Section \ref{sec:proofmain}, this would lead to a convergence of \eqref{defpartialprocess}, again under higher moment conditions.

This paper is structured as follows. In Section \ref{sec:results}, we state and discuss our main assumptions and state our results. The proofs can be found in Section \ref{sec:proofs} and a lengthy variance computation is contained in Section \ref{sec:calculations}. 
\\

\textbf{Acknowledgments:} The author is very grateful to Maurice Duits for several helpful discussions and for inspiring the author to investigate the partial traces.

\section{Random ensembles and main results}
\label{sec:results}

Let us begin with a closer look at the partial trace. When $Z=Z^\sn$ is a $n\times n$ complex Hermitian matrix, by the spectral theorem we may write $Z^\sn = U^\sn \Lambda^\sn (U^\sn)^*$, where $U^\sn$ is a $n\times n$ unitary matrix, $\Lambda^\sn$ is real diagonal with the eigenvalues $\lambda_1,\dots, \lambda_n$ on the diagonal and $A^*$ denotes the conjugate transpose of $A$. If $Z^\sn$ is real symmetric, $U^\sn$ is orthogonal instead. With this decomposition, we have for the partial trace as defined in \eqref{defpartialtrace},
\begin{align}\label{partialtrace2}
X_t^\sn(f) = \sum_{i=1}^{\lfloor tn\rfloor } f(Z)_{i,i} = \sum_{i=1}^{\lfloor tn\rfloor } \big(U^\sn f(\Lambda^\sn) (U^\sn)^*\big)_{i,i} 
= \sum_{i=1}^{\lfloor tn\rfloor } \sum_{k=1}^n |U^\sn_{i,k}|^2 f(\lambda_k) . 
\end{align}
The main object of our study is then the random non-negative finite measure $X_t^\sn$ defined by
\begin{align} \label{mainmeasure}
X_t^\sn =  \sum_{k=1}^n w_{k,t}^{(n)} \delta_{\lambda_k} . 
\end{align} 
with $\delta_z$ the Dirac measure in $z$ and the weights are given by 
\begin{align} \label{defweights}
w_{k,t}^{(n)} = \sum_{i=1}^{\lfloor tn\rfloor }  |U^\sn_{i,k}|^2 .
\end{align}
In this case $\mu(f)$ is just the shorthand notation for $\int f \, d\mu$. Note that the total mass of $X_t^\sn $ is given by $\lfloor tn\rfloor$. The representation \eqref{mainmeasure} shows that statements about the partial trace are in fact statements about a weighted version of the classical empirical eigenvalue distribution, which we denote by $\hat\mu^\sn$ and which corresponds to all weights being equal to $n^{-1}$. In \eqref{mainmeasure}, the weight of $\lambda_i$ is a norm of the first $\lfloor tn\rfloor$ entries of the corresponding eigenvector. Setting $t=1$, all weights in \eqref{mainmeasure} become 1, so that $\hat\mu^\sn = \frac{1}{n}X_1^\sn $. In other words, $n \hat\mu^\sn (f)$ is the linear eigenvalue statistic. 

Another prominent eigenvalue measure is the spectral measure $\mu_1^\sn$ of the pair $(Z^\sn ,e_1)$, defined by the functional calculus relation $\mu_1^\sn(f) = e_1^*f(Z^\sn)e_1 = f(Z^\sn)_{1,1}$. That is, the CLT in \eqref{CLTentry} is in fact a statement about $\mu_1^\sn(f)$. The spectral measure can be obtained from the partial trace by $\mu_1^\sn= X_{1/n}^\sn$. Although for classical ensembles of random matrices, the measures $\hat\mu^\sn$ and $\mu_1^\sn$ have the same limit in probability as $n\to \infty$, the fluctuations around this limit are very different, which becomes evident in the different central limit theorems in \eqref{CLTtrace} (for $n \hat\mu^\sn$) and in \eqref{CLTentry} (for $\mu_1^\sn$). The additional randomness of the weights in the spectral measure leads to substantially larger fluctuations. Let us remark that a similar behavior can be observed on the scale of large deviations: while $\hat\mu^\sn$ satisfies a large deviation principle with speed $n^2$ , see \cite{arous1997large} or \cite{anderson2010introduction}, for $\mu_1^\sn$ this is reduced to speed $n$ \cite{gamboa2011large,gamboa2016sum}.

\subsection{Assumptions}
\label{sec:results1}

In order to present the results for complex and real matrices in a unified expression, we follow the classical notation of \cite{dyson1962threefold} and introduce the parameter $\beta$, where $\beta=1$ if $U^\sn$ is real and orthogonal and $\beta=2$ if $U^\sn$ is complex and unitary. Let $\beta'=\beta/2$. 
We will always make the following assumption:
\begin{itemize}
\item[(A1)] \label{A1}The matrices $U^\sn$ and $\Lambda^\sn$ are independent and $U^\sn$ is Haar distributed on the unitary group ($\beta=2$) or on the orthogonal group ($\beta=1$).
\end{itemize}
Under assumption (A1), we can write the distribution of $(U^{(n)},\Lambda^{(n)})$ as $\mathbb{P} = \mathbb{P}_H \otimes \mathbb{P}_\Lambda$, where $\mathbb{P}_H$ is the Haar measure on the unitary group and $ \mathbb{P}_\Lambda$ is the distribution of the eigenvalues.  We denote expectation with respect to $\mathbb{P}_H$ and $\mathbb{P}_\Lambda$ by $\mathbb{E}_H$ and $\mathbb{E}_\Lambda$, respectively. Any convergence in distribution will be under $\mathbb{P}$ unless we specify otherwise. 
Without loss of generality, assume that all matrices $(U^{(n)},\Lambda^{(n)})$ for $n\geq 1$ are defined on a common probability space. While the distribution of $U^\sn$ is completely specified by (A1), we need that the empirical measure of the eigenvalues converges to a deterministic limit. Apart from Theorem \ref{thm:mainbridge}, we also assume a CLT for the linear eigenvalue statistic. Note that the two next assumptions are also conditions on the test function $f:\mathbb{R} \to \mathbb{R}$. 
\begin{enumerate}
\item[(A2)]  There exists a deterministic probability measure $\nu$, such that $\hat\mu^\sn$ converges weakly to $\nu$ $\mathbb{P}_\Lambda$-almost surely. Furthermore, $\hat\mu^\sn(f)$ converges to $\nu(f) $ and $\hat\mu^\sn(f^2)$ converges to $\nu(f^2)$.
\item[(A3)]  
There exists a $\sigma_1^2(f)\in[0,\infty)$ such that
\begin{align*}
X_1^\sn (f) - \mathbb{E} [X_1^\sn (f)] \xrightarrow[n \rightarrow \infty ]{d} \mathcal{N}(0,\sigma_1^2(f)) .
\end{align*}
\end{enumerate}

Let us comment on the assumptions above. Suppose the matrix $Z^\sn$ is distributed with density proportional to 
\begin{align} \label{densityV}
\exp \{ - \tfrac{1}{2}n\beta \tr V(X) \} 
\end{align}
with respect to the Lebesgue measure in each independent real entry in $X$. The potential $V:\mathbb{R} \to (-\infty,\infty]$ is supposed to be continuous and satisfy the growth (or confinement) condition
\begin{align} \label{confinement}
\liminf_{|x| \to \infty} \frac{V(x)}{2 \log |x|} > 1.  
\end{align}
The density \eqref{densityV} implies that assumption (A1) is satisfied and that the eigenvalues have a joint density proportional to 
\begin{align} \label{densityev}
\prod_{i<j} |\lambda_i-\lambda_j|^{\beta} \prod_{i=1}^n \exp\{ - \tfrac{1}{2}n\beta V(\lambda_i) \} 
\end{align}
with respect to the Lebesgue measure on $\mathbb{R}^n$ (see \cite{mehta2004random}). It follows from the large deviation principle of \cite{arous1997large} that the empirical eigenvalue distribution $\hat\mu^\sn$ converges exponentially fast to a compactly supported measure $\nu$.
Since the probability of deviating from the limit in the weak topology decays exponentially fast, the weak convergence holds almost surely on any joint probability space. That is, assumption (A2) is satisfied for any $f$ continuous and bounded. If moreover $\nu$ is supported by a single interval and the effective potential
\begin{align}\label{effpot}
\mathcal{J}_V(x) = \tfrac{1}{2}V(x)- \int \log | x-\xi| \, d\nu(\xi)
\end{align}
attains its infimum only on the support of $\nu$, then the largest and smallest eigenvalues each satisfy a large deviation principle \cite{arous2001aging,albeverio20011}. This implies that the probability of the extremal eigenvalues deviating from the support of $\nu$ decays exponentially and one easily obtains that (A2) holds also for continuous $f$ growing at most polynomially at infinity.

Turning to assumption (A3), it was shown in \cite{bekerman2018clt} that the CLT holds for quite general $\beta$-ensembles with density \eqref{densityev}, when $\mathcal{J}_V$ attains its infimum only on the support of $\nu$, and when $V$, $f$ are sufficiently smooth. We remark that \cite{bekerman2018clt} center by $\int f\, d\nu$, but their control of exponential moments allows to center as in (A3).  
For the classical cases of the Gaussian orthogonal ensemble (GOE, $\beta=1$) and Gaussian unitary ensemble (GUE, $\beta=2$) with potential $V(x) = \tfrac{1}{2}x^2$, the limit measure $\nu$ is the semicircle law with Lebesgue density $\tfrac{1}{2\pi}\sqrt{4-x^2}\mathbbm{1}_{[-2,2]}(x)$, and (A3) holds for any $f\in C^1(\mathbb{R})$ growing at most polynomially, with limiting variance
\begin{align}
\sigma_1^2(f) = \frac{1}{2\beta \pi^2} \int_{-2}^2 \int_{-2}^2 \left( \frac{f(x)-f(y)}{x-y}\right)^2 \frac{4-xy}{\sqrt{4-x^2} \sqrt{4-y^2} }\, d x d y ,
\end{align}
see \cite{pastur2011eigenvalue}, Chapter 3.2. 

As already mentioned in the introduction, the CLT in (A3) (and also assumption (A2)) is not only proven for random matrices with density \eqref{densityV}, but for a large variety of models, for example general Wigner or Wishart matrices. Such matrices have in general no Haar distributed matrix of eigenvectors, such that assumption (A1) fails to hold. However, given a random matrix $Z^\sn$ satisfying (A2) and (A3), we may take $U^\sn$ Haar distributed on the orthogonal or unitary group, and define $\tilde Z^\sn=U^\sn Z^\sn (U^\sn)^*$. Then the matrix $\tilde Z^\sn$ trivially has a Haar distributed matrix of eigenvectors independent of the eigenvalues. The second and third assumption continue to hold, so that now $\tilde Z^\sn$ satisfies all asumptions. 

Finally, let us remark that the method in the present paper also works if a weak convergence as in (A3) holds with a non-Gaussian limit, but to stay within the framework of CLTs for the linear eigenvalue statistic, we restrict the presentation to the Gaussian case.


\subsection{Results}

The following first theorem can be seen as a preview on the process convergence in Theorem \ref{thm:main} and highlights already the different effects the weights and eigenvalues have on the fluctuations. It shows a CLT for the weighted spectral measure or a single entry of the trace, more precisely, the asymptotic normality of  
\begin{align} \label{decomposition}
\sqrt{n}(\mu_{1}^\sn (f) - \mathbb{E}[ \mu_1^\sn (f)]) = \sqrt{n}(\mu_{1}^\sn (f) -  \hat \mu^\sn (f)) +\sqrt{n}( \hat \mu^\sn (f) - \mathbb{E}[\hat \mu^\sn (f)]) ,
\end{align}
where we recall that $\mu_{1}^\sn = X_{1/n}^\sn$ and $\hat\mu^\sn = \frac{1}{n}X_1^\sn $, as defined in the beginning of Section \ref{sec:results}. The random weights are responsible for the weak convergence of the first term on the right hand side, while under (A3) the second term has fluctuations of smaller order, and vanishes in the limit. Moreover, although both terms depend on the eigenvalues, they are asymptotically independent.

\begin{thm} \label{thm:spectralCLT}
Assume that (A1) and (A2) hold with $f\in \mathcal{C}_b(\mathbb{R})$, then
\begin{align*}
\sqrt{n}(\mu_{1}^\sn (f) - \hat \mu^\sn (f)) \xrightarrow[n \rightarrow \infty ]{d}  \mathcal{N}(0,\sigma_0^2(f)) , 
\end{align*}
where $\sigma_0^2(f) = \tfrac{2}{\beta}(\nu(f^2)-\nu(f)^2)$. 
If additionally $\sqrt{n} (\hat \mu^\sn (f)-\mathbb{E}[\hat \mu^\sn (f)])\xrightarrow[n \rightarrow \infty ]{d}  \mathcal{N}(0,\hat\sigma^2(f))$ with $\hat\sigma^2(f)\in [0,\infty)$, then
\begin{align*}
\sqrt{n}(\mu_{1}^\sn (f) - \mathbb{E}[\mu_1^\sn (f)]) \xrightarrow[n \rightarrow \infty ]{d}  \mathcal{N}(0,\sigma_0^2(f)+\hat\sigma^2(f)) . 
\end{align*}
In particular, if (A3) holds, then this convergence follows with $\hat\sigma^2(f)=0$. 
\end{thm}

Let us remark that for $Z$ a random matrix satisfying (A1) and (A2), the first convergence in Theorem \ref{thm:spectralCLT} may be rewritten as 
\begin{align} \label{firstentryconv1}
\sqrt{n}(f(Z)_{1,1} - \tfrac{1}{n} \tr f(Z)) \xrightarrow[n \rightarrow \infty ]{d}  \mathcal{N}(0,\sigma_0^2(f))
\end{align}
and if the distribution of $Z$ satisfies also (A3), then the second convergence in Theorem \ref{thm:spectralCLT} is equivalent to \eqref{CLTentry}.

As described in the introduction, the main objective is to show how the fluctuations of the linear eigenvalue statistic emerges from summing individual matrix elements. So now we consider the process
\begin{align}\label{defnormprocess}
\mathcal{X}^\sn(f) = \big( X_t^\sn(f)- \mathbb{E}[X_t^\sn(f)] \big)_{t\in [0,1]}
\end{align}
as a random element of $\mathcal{D}[0,1]$, equipped with the Skorokhod-topology and the Borel-$\sigma$ algebra. Our main result is then the following theorem.

\begin{thm}\label{thm:main}
Under assumptions (A1), (A2) and (A3), the process $\mathcal{X}^\sn(f) $ converges as $n\to \infty$ in distribution towards the continuous centered Gaussian process $\mathcal{X}(f)$ with covariance 
\begin{align*}
 \operatorname{Cov}(\mathcal{X}_s(f),\mathcal{X}_t(f)) = (t\wedge s -ts) \sigma_0^2 (f)  + ts \sigma_1^2(f) ,
\end{align*}
with $\sigma_0^2(f)$ as in Theorem \ref{thm:spectralCLT} and $\sigma_1^2(f)$ as in (A3). 
\end{thm}

The proof of Theorem \ref{thm:main} relies on a decomposition of the process $\mathcal{X}^\sn(f)$ into a sum of two processes similar to \eqref{decomposition}. We have 
\begin{align} \label{decomposition2}
\mathcal{X}^\sn(f) = \mathcal{W}^\sn(f) + \mathcal{Z}^\sn(f) ,
\end{align}
where 
\begin{align} \label{defWprocess}
\mathcal{W}^\sn_t(f) = \mathcal{X}_t^\sn(f) - \E_H \left[ \mathcal{X}_t^\sn(f)\right] , \qquad 0\leq t\leq 1 
\end{align}
is the process centered with respect to $\mathbb{P}_H$ and 
\begin{align} \label{defZprocess}
\mathcal{Z}^\sn_t(f) = \E_H \left[ \mathcal{X}_t^\sn(f)\right] - \E \left[  \mathcal{X}_t^\sn(f)\right] , \qquad 0\leq t\leq 1 .
\end{align}
Since $\E_H[|U_{i,k}|^2] = 1/n$, we have by \eqref{mainmeasure} and \eqref{defweights} 
\begin{align} \label{defWprocess2}
\mathcal{W}^\sn_t(f) = \sum_{k=1}^n \sum_{i=1}^{\lfloor tn\rfloor} \big( |U_{i,k}|^2-\tfrac{1}{n} \big) f(\lambda_k) 
\end{align} 
and  
\begin{align} \label{defZprocess2}
\mathcal{Z}^\sn_t(f) = \frac{\lfloor tn\rfloor}{n} \big( X_n(f) -\E [X_n(f)]\big) .
\end{align}
This decomposition has a similar effect as \eqref{decomposition}. The elements of the unitary matrix $U$ are the main source for the fluctuations of $\mathcal{W}^\sn(f)$ and this process is asymptotically independent of $\mathcal{Z}^\sn(f)$. Since by assumption (A3), $ \mathcal{Z}_t^\sn(f)$ converges to a Gaussian multiplied by $t$, this will result in the convergence of the sum. The main step in the proof of Theorem \ref{thm:main} is then the following functional limit theorem for the process $\mathcal{W}^\sn(f)$. Note that assumption (A3) is not needed for this part.

\begin{thm}\label{thm:mainbridge}
Suppose (A1) and (A2) are satisfied. Then $\mathbb{P}_\Lambda$ almost surely, as $n\to \infty$, the process $\mathcal{W}^\sn(f)$ converges in distribution under $\mathbb{P}_H$ towards $\sigma_0(f) B$, where $B$ is a standard Brownian bridge. 
\end{thm}

\begin{remark} \label{rem:quenched}
The weak convergence in Theorem \ref{thm:mainbridge} can be seen as a \emph{quenched} convergence, valid for almost all realizations of sequences of eigenvalues. It demonstrates that after centering with respect to $\P_H$, the origin of the random fluctuations of the partial trace is solely in the weights \eqref{defweights}, that is, in the eigenvector matrix. The eigenvalues give only a deterministic contribution in the limit, depending only on the equilibrium measure $\nu$. It is therefore not relevant for the result that the $\lambda_i$ are eigenvalues of a random matrix. Instead, Theorem \ref{thm:mainbridge} holds for any randomly weighted measure as in \eqref{mainmeasure} with weights \eqref{defweights}. For example, one could replace the support points of this measure with i{.}i{.}d{.} random variables, or realizations of a point process, as long as they are independent of the weights and assumption (A2) holds.
The same remark can be made about the first convergence in Theorem \ref{thm:spectralCLT}. It does not require (A3) and although it is not explicitly stated, the convergence holds under $\P_H$ for $\P_\Lambda$ almost all support points of the random measure $\mu_1^\sn$. 
\end{remark}

\section{Proofs}
\label{sec:proofs}

\subsection{Proof of Theorem \ref{thm:spectralCLT}}
\label{sec:proofbridge}

Let $\beta'=\beta/2$. It follows from the Haar distribution of $U^\sn$ that the vector of weights $(|U^\sn_{1,1}|^2, \dots ,|U^\sn_{1,n}|^2)$ has a  homogeneous Dirichlet distribution $\operatorname{Dir}_n(\beta')$, which is defined by the Lebesgue density for the first $n-1$ coordinates proportional to 
\begin{align*}
\big( x_1 \cdots x_{n-1}(1-x_1-\dots -x_{n-1}) \big)^{\beta' -1} \mathbbm{1}_{\{ x_i > 0, x_1+\dots +x_{n-1}<1 \} } . 
\end{align*}
The uniform distribution on the standard simplex corresponds thus to $\beta=2$. We will prove the CLT for weights following the general distribution $\operatorname{Dir}_n(\beta')$ for any $\beta'>0$, since it makes no difference in the proof. The starting point is the observation that the Dirichlet distribution can be generated by self-normalizing a vector of independent gamma random variables. More precisely, let $\gamma_1, \dots, \gamma_n$ be independent random variables with distribution $\operatorname{Gamma}(\beta')$, then  
\begin{align}\label{gammadir}
\left(\frac{\gamma_1}{\gamma_1+ \dots + \gamma_n}, \dots , \frac{\gamma_n}{\gamma_1+ \dots + \gamma_n}\right)
\sim \operatorname{Dir}_n(\beta'). 
\end{align}
where $\gamma_1, \dots, \gamma_n$ are independent and identically gamma distributed with parameters $(\beta',1)$ and mean $\beta'$, with moment generating function given for $t<1$ as
\begin{align}
\mathbb{E} \big[ e^{t \gamma_1} \big] = (1-t)^{-\beta'}.
\end{align}

Define the non-negative measure 
\begin{align} \label{defmutilde}
\tilde{\mu}_1^\sn = \frac{1}{n\beta'}\sum_{k=1}^n \gamma_k \delta_{\lambda_k} , 
\end{align}
then by \eqref{gammadir} the normalized measure $\tilde{\mu}_1^\sn \cdot \tilde{\mu}_1^\sn(1)^{-1}$ has the same distribution as 
${\mu}_1^\sn $. We first prove the convergence with $\mu_1^\sn$ replaced by $\tilde{\mu}_1^\sn$. Assume without loss of generality that $\nu(f)=0$. The moment generating function with respect to $\mathbb{P}_H$ is
\begin{align} \label{charfct}
& \mathbb{E}_H\left[ \exp\left\{ t \sqrt{n} (\tilde{\mu}_1^\sn(f)-\hat{\mu}^\sn(f) ) \right\} \right]
 =  \prod_{k=1}^n\mathbb{E}_H\left[ \exp\left\{  t (\sqrt{n}\beta')^{-1} \gamma_k f(\lambda_k) \right\} \right] \exp\{ -t \sqrt{n}^{-1}f(\lambda_k) \}  \notag  \\
& \qquad \qquad \qquad \qquad =  \prod_{k=1}^n \big( 1- t (\sqrt{n}\beta')^{-1} f(\lambda_k) \big)^{-\beta'} \exp\{ -t \sqrt{n}^{-1}f(\lambda_k) \}  \notag  \\
& \qquad \qquad \qquad \qquad = \exp \left\{ \sum_{k=1}^n  \left( -\beta' \log\big( 1- t (\sqrt{n}\beta')^{-1} f(\lambda_k) \big) -t \sqrt{n}^{-1}f(\lambda_k) \right)\right\} , 
\end{align}
where we used the independence of the weights and the independence of weights and eigenvalues and we take $|t|< \beta'||f||^{-1}_\infty$. Expanding the logarithm as $\log(1+x) = x -x^2/2 + r(x)$ with $|r(x)|\leq |x|^3 $ for $|x|\leq 1/2$ this gives
\begin{align} \label{charfct2}
& \mathbb{E}_H \left[ \exp\left\{ t \sqrt{n} (\tilde{\mu}_1^\sn(f)-\hat{\mu}^\sn(f) ) \right\} \right]
 = \exp \left\{ t^2/2 (\beta')^{-1} \hat{\mu}^\sn(f^2) + R_n(t,f) \right\} , 
\end{align}
with $|R_n(t,f)|\leq \sqrt{n}^{-1} (\beta')^2 |t|^3||f||_\infty^3$ for $n$ large enough. By Assumption (A2),  $\hat{\mu}^\sn(f^2)$ converges to $\nu(f^2)=\nu(f^2)-\nu(f)^2$ almost surely with respect to $\mathbb{P}_\Lambda$. Since $\hat{\mu}^\sn(f^2)$ and $R_n(t,f)$ are uniformly bounded (for $t$ and $f$ fixed), we have by dominated convergence
\begin{align}
\lim_{n\to \infty } \mathbb{E}\left[ \exp\left\{ t \sqrt{n} (\tilde{\mu}_1^\sn(f)-\hat{\mu}^\sn(f) ) \right\} \right]
= \exp \left\{ t^2/2 (\beta')^{-1} \nu(f^2) \right\} ,
\end{align}
that is, 
\begin{align} \label{spectralCLTtilde}
 \sqrt{n} (\tilde{\mu}_1^\sn(f)-\hat{\mu}^\sn(f) )  \xrightarrow[n \rightarrow \infty ]{d}  \mathcal{N}(0,(\beta')^{-1} \nu(f^2)) . 
\end{align}

In order to come back to the original measure $\mu_1^\sn(f) = \tilde{\mu}_1^\sn \cdot \tilde{\mu}_1^\sn(1)^{-1}$ we write
\begin{align} \label{spectralCLTtilde2}
 \sqrt{n} \big(\mu_1^\sn(f)-\hat{\mu}^\sn(f) \big) 
 =  \sqrt{n} \big(\tilde{\mu}_1^\sn(f)-\hat{\mu}^\sn(f) \big)\tilde{\mu}_1^\sn(1)^{-1} + \sqrt{n} \hat{\mu}^\sn(f)\big( \tilde{\mu}_1^\sn(1)^{-1} -1 \big) . 
\end{align}
By the strong law of large numbers, $\tilde{\mu}_1^\sn(1)$ converges almost surely to $\mathbb{E}[\tilde{\mu}_1^\sn(1)]= \mathbb{E}[(\beta')^{-1} \gamma_1] = 1$. So to conclude the convergence \eqref{spectralCLTtilde} with $\tilde{\mu}_1^\sn$ replaced by ${\mu}_1^\sn$, it suffices to show that the second term in \eqref{spectralCLTtilde2} vanishes in probability. Since $\hat{\mu}^\sn(f)$ converges almost surely to $\nu(f)=0$, this will follow if $\sqrt{n}(\tilde{\mu}_1^\sn(1)-1)$ is bounded in $L^2(\mathbb{P})$, which is easily checked by 
\begin{align}
n \mathbb{E}\big[ (\tilde{\mu}_1^\sn(1)-1)^2 \big] = n \mathbb{E} \left[ \left( \frac{1}{n} \sum_{i=1}^n ((\beta')^{-1}\gamma_i-1)\right)^2\right] 
 =  (\beta')^{-2} \operatorname{Var}(\gamma_1) = (\beta')^{-1} . 
\end{align}
This implies then that the last term in \eqref{spectralCLTtilde2} vanishes in probability and then by  \eqref{spectralCLTtilde}  the left hand side converges to $ \mathcal{N}(0,(\beta')^{-1} \nu(f^2))$ in distribution. This proves the first convergence in Theorem \ref{thm:spectralCLT}.

The second convergence in Theorem \ref{thm:spectralCLT} will follow from Lemma \ref{lem:composite} below. To apply it to the present setting, we may set $\mathbb{P}_1=\mathbb{P}_H$, $\mathbb{P}_2= \mathbb{P}_\Lambda$, 
\begin{align}
X^\sn = \sqrt{n} \big(\mu_1^\sn(f)-\hat{\mu}^\sn(f) \big), \qquad Y^\sn = \sqrt{n} (\hat \mu^\sn (f)-\mathbb{E}[\hat \mu^\sn (f)]) .
\end{align}
By assumption, $Y^\sn$ converges in distribution under $\mathbb{P}_\Lambda$ to $Y\sim \mathcal{N}(0,\hat\sigma^2(f))$. From \ref{charfct2} we get $\P_\Lambda$-almost surely
\begin{align}
\lim_{n\to \infty} \mathbb{E}_H\left[ \exp\left\{ t X^\sn \right\} \right] = \exp \left\{ t^2/2 (\beta')^{-1} \hat{\mu}^\sn(f^2)  \right\} ,
\end{align}
for any $t\in (-\beta'||f||^{-1}_\infty,\beta'||f||^{-1}_\infty )$. Since the moment generating functions are continuous, almost sure convergence for fixed $t$ implies almost sure pointwise convergence, which implies that the convergence \eqref{conditionalconv} holds with $X\sim \mathcal{N}(0,\sigma_0^2(f))$. 
Lemma \ref{lem:composite} implies then the convergence of $X^\sn+Y^\sn$ to $X+Y$. Noting that $\mathbb{E}[\hat \mu^\sn (f)]=\mathbb{E}[ \mu_1^\sn (f)]$, this finishes the proof. 
\QED

\begin{lem} \label{lem:composite}
Let $(\Omega_1\times\Omega_2, \mathcal{G}, \mathbb{P}_1\otimes \mathbb{P}_2)$ be a probability space and $X^\sn:\Omega_1\times\Omega_2\to \Omega'$ and $Y^\sn:\Omega_2 \to \Omega'$ random variables, where $\Omega'$ is a separable metric space with Borel $\sigma$-algebra. If $Y^\sn$ converges to $Y$ in distribution under $\P_2$ and 
\begin{align}\label{conditionalconv}
\mathbb{E}_1[F(X^\sn)] \xrightarrow[n \rightarrow \infty ]{} \mathbb{E}[F(X)]
\end{align}
$\mathbb{P}_2$-almost surely for any bounded continuous $F:\Omega'\to \mathbb{R}$, where $\mathbb{E}_1,\E$ is the expectation with respect to $\mathbb{P}_1, \P_1\otimes \P_2$ respectively, then
\begin{align}
(X^\sn,Y^\sn) \xrightarrow[n \rightarrow \infty ]{d} (X,Y)
\end{align}
in distribution under $\P_1\otimes \P_2$, with $X$ and $Y$ independent.
\end{lem}

\textbf{Proof:} The main observation is that functions $\mathcal{F}:\Omega'\times \Omega'\to \mathbb{R}$ with $\mathcal{F}(x,y)=F(x)G(y)$ and $F,G$ bounded continuous are sufficient to determine convergence in distribution, see Lemma 4.1 in \cite{hoffmann1977probability}. For such $F,G$, we have
\begin{align*}
& \ \quad \big| \E[F(X^\sn)G(Y^\sn)] - \E[F(X)]\E[G(Y)] \big| \\
& \leq \big| \E[(F(X^\sn)-\E[F(X)]) G(Y^\sn)]\big|  + \big| \E[F(X)] (\E[G(Y^\sn)]  - \E[G(Y)]) \big| \\
& = \big| \E_2[(\E_1[F(X^\sn)]-\E[F(X)]) G(Y^\sn)]\big| +  \big| \E[F(X)] (\E[G(Y^\sn)]  - \E[G(Y)]) \big| .
\end{align*}
The first term vanishes by dominated convergence using \eqref{conditionalconv}, the second one by the convergence of $Y^\sn$ under $\P_2$. 
\QED

\subsection{Proof of Theorem \ref{thm:mainbridge}}
\label{sec:proofbridge}

\subsubsection{Representation by a bivariate process}

We will first show the statement of Theorem \ref{thm:mainbridge} for piecewise constant functions $h$ with
\begin{align}\label{elementaryh}
h(\lambda) = \sum_{m=1}^M \gamma_m \mathbbm{1}_{ (a_{m},b_{m}]}(\lambda),
\end{align}
for some real $\gamma_m, 1\leq m\leq M$, and $a_1<b_1\leq a_2< \dots \leq b_m$ such that $\nu((-\infty,\cdot])$ is continuous at all $a_i,b_i$. 
The last condition only excludes countable many points for the choice of $a_i,b_i$ and in particular still allows to approximate any $f\in L^2(\nu)$. 
Let $U^{(n)}$ be a sequence of $n\times n$ unitary or orthogonal Haar distributed matrices. We denote by $\widetilde{\mathcal{W}}^\sn$ a process indexed by subsets $A\times B$ of $\{1,\dots ,n\}^2$, such that
\begin{align*}
\widetilde{\mathcal{W}}^\sn_{A,B} = \sum_{i,j=1}^n \big( |U^{(n)}_{i,j}|^2 - \tfrac{1}{n}\big) \mathbbm{1}_A(i)\mathbbm{1}_B(j) .
\end{align*}
If $A$ and/or $B$ are of the form $\{1,\dots ,\lfloor tn\rfloor \}$ with $t\in [0,1]$, we replace the corresponding index by $t$.  

We consider $(\widetilde{\mathcal{W}}^\sn_{s,t})_{s,t\in[0,1]}$ as a random element of $\mathcal{D}([0,1]^2)$, the multidimensional version of the Skorokhod-space. $\mathcal{D}([0,1]^2)$ contains all $X:[0,1]^2\to \mathbb{R}$ which are ``continuous from the north-east'' and have existing limits in each quadrant, i.e., $\lim_{s\searrow s_0,t\searrow t_0} X(s,t)= X(s_0,t_0)$ and $\lim_{s\searrow s_0,t\nearrow t_0} X(s,t)$, $\lim_{s\nearrow s_0,t\searrow t_0} X(s,t)$ and $\lim_{s\nearrow s_0,t\nearrow t_0} X(s,t)$ exist. 
We endow $\mathcal{D}([0,1]^2)$ with a generalization of Skorokhod's $J_1$-metric defined by
\begin{align} \label{skorokhodmetric}
d(X,Y) = \inf_{\lambda_1,\lambda_2} \max \left\{ \sup_{s\in[0,1]} |\lambda_1(s)-s| ,\sup_{t\in[0,1]} |\lambda_2(t)-t| , \sup_{s,t\in[0,1]} |X(\lambda_1(s),\lambda_2(t))-Y(s,t)|\right\} , 
\end{align}
where the infimum is taken over all continuous one-to-one mappings $\lambda_i:[0,1]\to [0,1]$ fixing $0$. Then as in the one-dimensional case, $\mathcal{D}([0,1]^2)$ with metric \eqref{skorokhodmetric} is separable and although it is not complete, there is an equivalent metric such that $\mathcal{D}([0,1]^2)$ becomes complete, see Section 5 of \cite{straf1972weak} or \cite{bickel1971convergence}, Section 3. 
As in the classical case layed out in Section 12 of \cite{billingsley1999convergence}, convergence with respect to the metric \eqref{skorokhodmetric} with a continuous limit actually implies convergence in supremum norm.

It was shown in \cite{donati2012truncations}, that for suitable index sets, $\widetilde{\mathcal{W}}^\sn$ converges to $\sqrt{2/\beta} \mathcal{B}$, where $\mathcal{B}$ is a bivariate tied-down Brownian bridge, a centered Gaussian process on $[0,1]^2$ with continuous paths and covariance 
\begin{align} \label{covarianceBB}
\mathbb{E}[\mathcal{B}(s,t)\mathcal{B}(s',t')]= (s\wedge s'-ss')(t\wedge t'-tt') .
\end{align}

\begin{thm}[\cite{donati2012truncations}, Thm 1.1] \label{DRconvergence}
As $n\to \infty$, the process  $(\widetilde{\mathcal{W}}^\sn_{s,t})_{s,t\in[0,1]}$ 
converges in distribution under $\mathbb{P}_H$ to $\sqrt{2/\beta} \mathcal{B}$, with $\mathcal{B}$ a bivariate tied down Brownian bridge. 
\end{thm}

Now consider $h$ as in \eqref{elementaryh}, then
\begin{align}\label{bivariate1}
\mathcal{W}^\sn_t(h) & = \sum_{k=1}^n \sum_{i=1}^{\lfloor tn\rfloor} \big( |U_{i,k}|^2-\tfrac{1}{n} \big) h(\lambda_k) \notag \\
& = \sum_{m=1}^M \gamma_m \sum_{k=1}^n \sum_{i=1}^n \big( |U_{i,k}|^2-\tfrac{1}{n} \big) \mathbbm{1}_{ (a_m,b_m]}(\lambda_k) \mathbbm{1}_{\{1,\dots ,\lfloor tn\rfloor  \} }(i) \notag \\
& =\sum_{m=1}^M \gamma_m \widetilde{\mathcal{W}}^\sn_{t,A_m} , 
\end{align}
with $A_m = \{ k\, |\, \lambda_k \in (a_m,b_m] \}$. 
For $a\in \mathbb{R}$, let 
\begin{align}\label{empriricallambda}
F^\sn(s) = \tfrac{1}{n}|\{\lambda_k^{(n)}: \lambda_k^{(n)}\leq a\}|
\end{align}
be the normalized number of eigenvalues $\leq s$, then we claim that  
\begin{align}\label{equalitybivariate}
\mathcal{W}^\sn(h) \stackrel{\mathbb{P}_H}{=} \widetilde{\mathcal{W}}^\sn(h)  := \sum_{m=1}^M \gamma_m \big( \widetilde{\mathcal{W}}^\sn_{\cdot, F^\sn(b_{m})} -  
\widetilde{\mathcal{W}}^\sn_{\cdot, F^\sn(a_{m})} \big)\, ,
\end{align}
where $\stackrel{\mathbb{P}_H}{=}$ denotes equality in distribution under $\mathbb{P}_H$. To see this, let $\pi$ by a permutation of $\{1,\dots ,n\}$ such that $\lambda_{\pi(1)}\leq \dots \leq \lambda_{\pi(n)}$. If $\Pi$ is the permutation matrix with entries $\Pi_{i,j} = \mathbbm{1}_{\pi(i)=j}$, then $\Pi$ is orthogonal. By the invariance of the Haar measure, we have $U\stackrel{\mathbb{P}_H}{=}U\Pi$, which implies that
\begin{align}\label{bivariate3}
\mathcal{W}_t^\sn(h) & \stackrel{\mathbb{P}_H}{=} \sum_{j=1}^n \sum_{i=1}^{\lfloor t n\rfloor} \big( |(U\Pi)_{i,j}|^2-\tfrac{1}{n} \big) h(\lambda_j) \notag \\
& = \sum_{j=1}^n \sum_{i=1}^{\lfloor t n\rfloor} \big( |U_{i,\pi^{-1}(j)}|^2-\tfrac{1}{n} \big) h(\lambda_j) \notag \\
& = \sum_{j=1}^n \sum_{i=1}^{\lfloor t n\rfloor} \big( |U_{i,j}|^2-\tfrac{1}{n} \big) h(\lambda_{\pi(j)}) \notag \\
& = \sum_{m=1}^M \gamma_m \sum_{j=1}^n \sum_{i=1}^{\lfloor t n\rfloor} \big( |U_{i,j}|^2-\tfrac{1}{n} \big) \mathbbm{1}_{\{a_m<\lambda_{\pi(j)}\leq b_m \}} ,
\end{align}
and the last line equals the right hand side of \eqref{equalitybivariate}. The equality in distribution in \eqref{bivariate3} holds also when both sides are viewed as a function of $t$, which implies \eqref{equalitybivariate}. 
We are now almost in the situation to apply Theorem \ref{DRconvergence}.

\subsubsection{A subordination argument}

Assumption (A2) implies the $\mathbb{P}_\Lambda$-almost sure convergence of $F^\sn(s)$ defined in \eqref{empriricallambda} to $F(s) = \nu((-\infty,s])$ for all $s\in S=\{a_1,b_1,\dots ,a_M,b_M\}$.  
Together with Theorem \ref{DRconvergence} this will yield the convergence of $\widetilde{\mathcal{W}}^\sn$ at random time points given by $F^\sn$, and we show in this section the convergence 
\begin{align}\label{subordinatedconv}
({\mathcal{W}}_t^\sn(h))_{t\in[0,1]}  \xrightarrow[n \rightarrow \infty ]{d} \big( \mathcal{W}_t(h)\big)_{t\in [0,1]} := \left( \sum_{m=1}^M \gamma_m \sqrt{\tfrac{2}{\beta}}\big( \mathcal{B}_{t,F(b_{m})} - \mathcal{B}_{t,F(a_{m})} \big) \right)_{t\in[0,1]} 
\end{align}
$\mathbb{P}_\Lambda$-almost surely in distribution under $\mathbb{P}_H$. Recall that by \eqref{equalitybivariate} we have ${\mathcal{W}}_t^\sn(h)  \stackrel{\mathbb{P}_H}{=} \widetilde{\mathcal{W}}_t^\sn(h)$. We defined all unitary $U^\sn, n\geq 1$, and therefore also all $\widetilde{\mathcal{W}}^\sn, n\geq 1$ on a common probability space. By the Skorokhod representation theorem, there exists a modification of this space, such that $\widetilde{\mathcal{W}}^\sn \to \sqrt{2/\beta} \mathcal{B}$ almost surely, with respect to a measure we again denote by $\mathbb{P}_H$. The product structure implied by assumption (A1) allows us to extend this to a product space with law $\mathbb{P}_H\otimes \mathbb{P}_\Lambda$ such that 
\begin{align}\label{bivariatejointly}
\big((\widetilde{\mathcal{W}}_{s,t}^\sn)_{s,t\in [0,1]} , (F^\sn(s))_{s\in S} \big)  \xrightarrow[n \rightarrow \infty ]{} 
\left( \left(\sqrt{\tfrac{2}{\beta}}\mathcal{B}_{s,t}\right)_{s,t\in [0,1]}, (F(s))_{s\in S} \right) 
\end{align}  
$\mathbb{P}_H\otimes \mathbb{P}_\Lambda$-almost surely in $\mathcal{D}([0,1]^2)\times \mathbb{R}^{2M}$. 
By \eqref{equalitybivariate}, we need to consider
\begin{align}\label{subordinatedconv2}
& \sup_{t\in[0,1]} \big| \widetilde{\mathcal{W}}_t^\sn(h)  - \mathcal{W}_t(h)\big| \notag \\
& =
 \sup_{t\in[0,1]} \left| \sum_{m=1}^M \gamma_m \big( \widetilde{\mathcal{W}}^\sn_{t,F^\sn(b_{m})} -  
\widetilde{\mathcal{W}}^\sn_{t,F^\sn(a_{m})} \big) - \sum_{m=1}^M \gamma_m \sqrt{\tfrac{2}{\beta}} \big( \mathcal{B}_{t,F(b_{m})} - \mathcal{B}_{t,F(a_{m})} \big)\right| \notag  \\
& \leq \sum_{m=1}^M |\gamma_m| \left(  \sup_{t\in[0,1]} \left|  \widetilde{\mathcal{W}}^\sn_{t,F^\sn(b_{m})} - \sqrt{\tfrac{2}{\beta}}\mathcal{B}_{t,F(b_{m})} \right| + \sup_{t\in[0,1]} \left| \widetilde{\mathcal{W}}^\sn_{t,F^\sn(a_{m})} - \sqrt{\tfrac{2}{\beta}}\mathcal{B}_{t,F(a_{m})} \right| \right) .
\end{align}
An individual supremum in \eqref{subordinatedconv2} can then be bounded as 
\begin{align}\label{subordinatedconv3}
& \sup_{t\in[0,1]} \left|  \widetilde{\mathcal{W}}^\sn_{t,F^\sn(s)} - \sqrt{\tfrac{2}{\beta}} \mathcal{B}_{t,F(s)} \right| \notag \\
& \leq \sup_{t\in[0,1]} \left|  \widetilde{\mathcal{W}}^\sn_{t,F^\sn(s)} - \sqrt{\tfrac{2}{\beta}}\mathcal{B}_{t,F^\sn(s)} \right| 
+\sqrt{\tfrac{2}{\beta}} \sup_{t\in[0,1]} \left|  \mathcal{B}_{t,F^\sn(s)} - \mathcal{B}_{t,F(s)} \right| 
\end{align} 
with $s\in S$. Since $\mathcal{B}$ is uniformly continuous, the convergence of Theorem \ref{DRconvergence} holds with respect to the supremum norm on $\mathcal{D}([0,1]^2)$, which implies that the first term in \eqref{subordinatedconv3} vanishes as $n\to \infty$. Since $F^\sn(s) \to F(s)$ for $s\in S$ and using again the uniform continuity of $\mathcal{B}$, the second term vanishes as well. By the bound in \eqref{subordinatedconv2} the convergence $\widetilde{\mathcal{W}}^\sn(h) \to \mathcal{W}(h)$ follows $\mathbb{P}_H\otimes \mathbb{P}_\Lambda$-almost surely in $\mathcal{D}([0,1])$. 
The product structure of the extended probability space implies then that for any bounded continuous $G$ we get
\begin{align}\label{subordinateconv4}
\mathbb{E}_H\big[G( \mathcal{W}^\sn(h)) \big] =\mathbb{E}_H\big[G( \widetilde{\mathcal{W}}^\sn(h)) \big] 
 \xrightarrow[n \rightarrow \infty ]{} \mathbb{E}\big[G( \mathcal{W}(h)) \big] 
\end{align}
$\mathbb{P}_\Lambda$-almost surely, that is, \eqref{subordinatedconv} holds.

Since $\mathcal{B}$ is a centered Gaussian process with continuous paths, the same holds for $\mathcal{W}(h)$. To calculate the covariance we first note that according to \eqref{covarianceBB},
\begin{align} \label{bivariate4}
& \notag \quad \mathrm{Cov}\left( \gamma_m \sqrt{\tfrac{2}{\beta}}(\mathcal{B}_{s,F(b_m) } -  
\mathcal{B}_{s,F(a_{m}) }),\gamma_\ell \sqrt{\tfrac{2}{\beta}}(\mathcal{B}_{t,F(b_\ell)}  -  
\mathcal{B}_{t,F(a_{\ell}) })\right) \\
& = \gamma_m\gamma_\ell \tfrac{2}{\beta} (s\wedge t-st)\big[ F(b_m)\wedge F(b_\ell) - F(b_m)F(b_\ell) -(F(b_m)\wedge F(a_\ell) - F(b_m)F(a_\ell)) \notag \\
& \qquad \qquad \qquad \qquad -(F(a_m)\wedge F(b_\ell) - F(a_m)F(b_\ell))+F(a_m)\wedge F(a_\ell) - F(a_m)F(a_\ell)\big] .
\end{align}
For $m\neq \ell$ the minima in \eqref{bivariate4} all cancel and this reduces to
\begin{align*}
& \quad \gamma_m\gamma_\ell \tfrac{2}{\beta}(s\wedge t-st)\big[ - (F(b_m)-F(a_m))(F(b_\ell)-F(a_\ell))\big] \\
& = - (s\wedge t-st) \tfrac{2}{\beta} \int \gamma_m \mathbbm{1}_{(a_m,b_m]} d\nu \cdot \int \gamma_\ell \mathbbm{1}_{(a_\ell,b_\ell]} d\nu,
\end{align*}
while for $m=\ell$ we get
\begin{align*}
& \quad \gamma_m^2 \tfrac{2}{\beta} (s\wedge t-st)\big[ (F(b_m)-F(a_m))- (F(b_m)-F(a_m))^2\big] \\
& =  (s\wedge t-st) \tfrac{2}{\beta} \left[ \int  \gamma_m^2 \mathbbm{1}_{(a_m,b_m]} d\nu - \left( \int \gamma_m \mathbbm{1}_{(a_m,b_m]} d\nu\right)^2 \right]. 
\end{align*}
Summing over $m,\ell$, this yields for the covariance
\begin{align} \label{bivariate5}
\mathrm{Cov} (\mathcal{W}_s(h),\mathcal{W}_t(h)) = (s\wedge t-st) \tfrac{2}{\beta} \left[ \int h^2d\nu - \left(\int hd\nu \right)^2 \right] .
\end{align}
That is, $\mathcal{W}(h)= \sigma_1 (h) B$, with $B$ a standard Brownian bridge. It remains to replace the elementary function $h$ as in \eqref{elementaryh} by an arbitrary $f$.

\subsubsection{Extension to general $f$}

Let $f\in L^2(\nu)$, $G$ be a bounded uniformly continuous functional from $\mathcal{D}([0,1])$ to $\mathbb{R}$, and $\varepsilon>0$. Denoting now by $d$ the Skorokhod $J_1$-metric on $\mathcal{D}([0,1])$, let $\delta <1$ be so small that $d(X,Y) \leq \delta$ implies $|G(X) - G(Y)| \leq \varepsilon$. 
In order to extend the convergence of $\mathcal{W}^\sn(h)$ with $h$ as in the previous sections replaced by $f$, we need an a-priori estimate on the distance of the processes $\mathcal{W}^\sn(h)$ and $\mathcal{W}^\sn(f)$. The proof is postponed to the end of this section.

\begin{lem} \label{distanceestimate} 
There exists a constant $c>0$, such that for $\eta>0$ and $g:\mathbb{R}\to \mathbb{R}$ measurable, 
\begin{align*}
\limsup_{n\to \infty} \mathbb{P}_H \left( \sup_{t\in [0,1]}|\mathcal{W}^\sn_t(g)|  > \eta\right) \leq \limsup_{n\to \infty}\frac{c}{\eta^2} (\hat \mu^\sn(g^2) - \hat\mu^\sn(g)^2) , 
\end{align*}
$\mathbb{P}_\Lambda$-almost surely. 
\end{lem}

Note that for $g$ satisfying (A2), the upper bound in Lemma \ref{distanceestimate} is equal to $c\sigma_0^2(g)/\eta^2$. 
We may approximate $f$ by a piecewise constant function $h=h_\varepsilon$ as in \eqref{elementaryh}, such that $||f-h||_{L^2(\nu)}\leq \delta^2\varepsilon$. We want to apply Lemma \ref{distanceestimate} with $g=f-h$, however in order to control the upper bound we need to control $\hat \mu^\sn(fh)$. For this, we write $f=f_+-f_-$ with $f_\pm \geq 0$, and assume the positive and negative part $f_+$ and $f_-$ is approximated by $h_+$ and $h_-$ respectively, with $h_\pm \geq 0$ and such that $h_\pm \leq f_\pm$. Then we can estimate
\begin{align} \label{distancepm}
\mu^\sn(fh) = \mu^\sn(f_+h_+)+\mu^\sn(f_-h_-)\geq \mu^\sn(h_+^2)+\mu^\sn(h_-^2) = \mu^\sn (h^2), 
\end{align}
such that
\begin{align} \label{distancepm2}
\mu^\sn((f-h)^2) 
 \leq \mu^\sn(f^2-h^2) .
\end{align}
By assumption (A2), $\mu^\sn(f^2)\to \nu(f^2)$ $\P_\Lambda$ almost surely and the elementary form of $h$ as in \eqref{elementaryh} implies $\mu^\sn(h^2)\to \nu(h^2)$ as well. This implies that \eqref{distancepm2} converges $\P_\Lambda$ almost surely to 
$\nu(f^2-h^2)\leq 2 ||f-h||_{L^2(\nu)}||f||_{L^2(\nu)}$. We have $\P_\Lambda$ almost surely
\begin{align} \label{distancepm3} 
\mathbb{E}_H [|G(\mathcal{W}^\sn(f))-G(\mathcal{W}^\sn(h_\varepsilon))|] & \leq \varepsilon + 2 ||G||_\infty \mathbb{P}_H\big( d(\mathcal{W}^\sn(f),\mathcal{W}^\sn(h_\varepsilon)) >\delta \big) \notag \\
& \leq \varepsilon + 2 ||G||_\infty \mathbb{P}_H\big( ||\mathcal{W}^\sn(f)-\mathcal{W}^\sn(h_\varepsilon)||_\infty >\delta \big) ,
\end{align}
so that we obtain from Lemma \ref{distanceestimate} with $g=f-h$ and \eqref{distancepm2} 
\begin{align} \label{weakconv1}
\limsup_{n\to \infty} \big| \mathbb{E}_H[G(\mathcal{W}^\sn(f))]-\mathbb{E}[G(\mathcal{W}(h_\varepsilon))] \big| & \leq \varepsilon +2||G||_\infty \delta^{-2} \nu(f^2-h^2) \notag \\
& \leq  \varepsilon +4||G||_\infty \delta^{-2} ||f-h||_{L^2(\nu)}||f||_{L^2(\nu)} \notag \\
& \leq \varepsilon +4\varepsilon ||G||_\infty ||f||_{L^2(\nu)}  .
\end{align}
Furthermore, if we set $\mathcal{W}(f) = \sigma_1 (f) B$, then $\mathcal{W}(h_\varepsilon)$ and $\mathcal{W}(f)$ are Gaussian processes with covariance $(s\wedge t -st)\sigma_1^2(h_\varepsilon)$ and $(s\wedge t -st)\sigma_1^2(f)$, respectively, and if $\varepsilon \to 0$ and then $h=h_\varepsilon\to f$ in $L^2(\nu)$, 
\begin{align}\label{weakconv2}
\mathcal{W}(h_\varepsilon) = \sigma_1 (h_\varepsilon) B \xrightarrow[\varepsilon \rightarrow 0 ]{} \sigma_1 (f) B=\mathcal{W}(f) .
\end{align}
The combination of \eqref{weakconv1} and \eqref{weakconv2} shows that we may replace $h$ in \eqref{subordinateconv4} by any $f\in L^2(\nu)$, so that $\mathcal{W}^\sn(f)$ converges to $\mathcal{W}(f)$ in distribution under $\mathbb{P}_H$, for $\mathbb{P}_\Lambda$-almost all $\lambda$. \QED

\textbf{Proof of Lemma \ref{distanceestimate}:} We write
\begin{align*}
\mathcal{W}_t^\sn(g) = \sum_{j=1}^n \sum_{i=1}^{\lfloor tn\rfloor} \big( |U_{i,j}|^2-\tfrac{1}{n} \big) g(\lambda_j)  = \sum_{i=1}^{\lfloor tn\rfloor} Y_{i,n}
\end{align*}
where 
\begin{align*}
Y_{i,n}= \sum_{j=1}^n \big( |U_{i,j}|^2-\tfrac{1}{n} \big) g(\lambda_j) .
\end{align*}
By the invariance of the Haar distribution, the vector of increments $(Y_{1,n},\dots ,Y_{n,n})$ is exchangable under $\mathbb{P}_H$, meaning that any permutation of the $Y_{i,n}$ has the same distribution. Corollary 2 in \cite{pruss1998maximal} shows that there exists a universal constant $c>0$, such that 
\begin{align*}
\mathbb{P}_H \left( \sup_{1\leq k\leq n} \left| \sum_{i=1}^{k} Y_{i,n}\right| > \eta \right) \leq c \mathbb{P}_H \left(  \left|\sum_{i=1}^{\lfloor n/2 \rfloor} Y_{i,n}\right| > \eta/c \right)
\end{align*}
and the right hand side can be bounded by $c^3\mathbb{E}_H[\mathcal{W}^\sn_{1/2}(g)^2]/\eta^2$. The calculations in Section \ref{sec:calculations}, more precisely taking the $\limsup$ in \eqref{varconvergence2}, show that this upper bound implies the statment of Lemma \ref{distanceestimate}.  \QED

\subsection{Proof of Theorem \ref{thm:main}}
\label{sec:proofmain}

After completing the proof of Theorem \ref{thm:mainbridge}, this proof follows from Lemma \ref{lem:composite}, as in the proof of Theorem \ref{thm:spectralCLT}. Set $\P_1=\P_H, \P_2=\P_\Lambda$, and $X^\sn = (\mathcal{W}_t^\sn(f))_{t\in [0,1]}$, $Y^\sn = (\mathcal{Z}_t^\sn(f))_{t\in [0,1]}$. Then by Theorem \ref{thm:mainbridge}, the convergence \eqref{conditionalconv} holds with limit $X = (\mathcal{W}_t(f))_{t\in [0,1]}$ and by assumption (A3), $Y^\sn$ converges in distribution under $\P_2$ to $Y=(t\mathcal{Z}(f))_{t\in [0,1]}$, with $\mathcal{Z}(f) \sim \mathcal{N} (0,\sigma_1^2(f))$ (recall \eqref{defZprocess2}). Then Lemma \ref{lem:composite} implies the convergence 
\begin{align}
(\mathcal{X}_t^\sn(f))_{t\in [0,1]} = \big( \mathcal{W}_t^\sn(f) + \mathcal{Z}_t^\sn(f) \big)_{t\in [0,1]}
\xrightarrow[n \rightarrow \infty ]{d}
\big( \mathcal{W}_t(f) + t\mathcal{Z}(f) \big)_{t\in [0,1]}
\end{align} 
under $\P$, with $\mathcal{W}_t(f)$ and $\mathcal{Z}(f)$ independent. This is the convergence claimed in Theorem \ref{thm:main}. 
\QED

\section{Calculation of the covariance}\label{sec:calculations}

In this section we prove that in the setting of Theorem \ref{thm:mainbridge} for $s,t \in [0,1]$,
\begin{align}\label{varconvergence}
\lim_{n\to \infty} \operatorname{Cov}_H \big(\mathcal{X}_s(f),\mathcal{X}_t(f) \big) =  (s\wedge t -st ) \sigma_0^2(f) 
\end{align}
$\P_\Lambda$-almost surely, where $\operatorname{Cov}_H$ denotes the covariance with respect to $\P_H$. This requires to compute some mixed moments of entries of the eigenvector matrix $U^\sn$, where for the sake of a lighter notation, we drop the superscript. 
We recall that if $U=(U_{i,j})_{i,j}$ is Haar distributed on the unitary ($\beta=2$) or the orthogonal ($\beta=1$) group, $(|U_{i,1}|^2,\dots ,|U_{i,n}|^2)$ is $\operatorname{Dir}_n(\beta')$ distributed. Each $|U_{i,j}|^2$ follows then a beta distribution with parameter $(\beta',\beta'(n-1))$ and  therefore
\begin{align} \label{Umoments1}
\E_H[|U_{i,j}|^2] = \frac{1}{n},\quad \E_H[|U_{i,j}|^4] = \frac{1+\beta'}{n(\beta'n+1)} . 
\end{align}
Moreover, if $j\neq k$, then $(U_{i,j},U_{i,k})$ is Dirichlet distributed with parameter $(\beta',\beta',\beta'(n-2))$, which implies
\begin{align} \label{Umoments2}
\E_H[|U_{i,j}|^2|U_{i,k}|^2] = \frac{\beta'}{n(\beta'n+1)} .  
\end{align}
If additionally $m\neq i$, then using 
\begin{align} \label{Umoments3}
\E_H[|U_{i,j}|^2]  = \sum_{m'=1}^n \E_H[|U_{i,j}|^2|U_{m',k}|^2] =  \E_H[|U_{i,j}|^2|U_{i,k}|^2] + (n-1)  \E_H[|U_{i,j}|^2|U_{m,k}|^2] ,
\end{align}
we see that by \eqref{Umoments1} and \eqref{Umoments2}
\begin{align} \label{Umoments4}
\E_H[|U_{i,j}|^2|U_{m,k}|^2] = \frac{(n-1)\beta'+1}{n(n-1)(\beta'n+1)} .
\end{align}
These identities can also be obtained from \cite{collins2006integration} or for $\beta=2$ from Proposition 4.2.3 of \cite{hiai2006semicircle}. Now let $s,t\in [0,1]$ and set $s_n = \lfloor sn \rfloor$, $t_n = \lfloor tn \rfloor$. If $s\in \{0,1\}$ (or $t\in \{0,1\}$), then $\mathcal{X}_s(f)=0$ (or $\mathcal{X}_t(f)=0$) and \eqref{varconvergence} is trivially true. So assume that $s,t \in(0,1)$, and $n$ is so large that $s_n,t_n\geq 1$. Without loss of generality , let $s \leq r$. Then we get by \eqref{mainmeasure} for the mixed moment with respect to $\mathbb{P}_H$  
\begin{align*}
\E_H \left[ X^\sn_s(f) X^\sn_t(f)\right] & = \E_H \left[ \left( \sum_{i=1}^n \sum_{l=1}^{s_n} |U_{l,i}|^2 f(\lambda_i) \right) \left( \sum_{j=1}^n \sum_{m=1}^{t_n} |U_{m,j}|^2 f(\lambda_j) \right)\right] \\
 & = \sum_{l=1}^{s_n} \sum_{m=1}^{t_n} \sum_{i,j=1}^n \E_H \left[ |U_{l,i}|^2 |U_{m,j}|^2\right]  f(\lambda_i)f(\lambda_j) . 
\end{align*}
This sum can be decomposed according to whether $l=m$ or not as
\begin{align} 
&\quad \sum_{l,m=1, l\neq m}^{s_n} \sum_{i,j=1}^n \E_H \left[ |U_{l,i}|^2 |U_{m,j}|^2\right]  f(\lambda_i)f(\lambda_j)  \label{cov1}\\
&+  \sum_{l=1}^{s_n} \sum_{i,j=1}^n \E_H \left[ |U_{l,i}|^2 |U_{l,j}|^2\right]  f(\lambda_i)f(\lambda_j)  \label{cov2} \\
&+  \sum_{l=1}^{s_n} \sum_{m=s_n+1}^{t_n}\sum_{i,j=1}^n \E_H \left[ |U_{l,i}|^2 |U_{m,j}|^2\right]  f(\lambda_i)f(\lambda_j)  . \label{cov3}
\end{align}
The inner sum in \eqref{cov1} and \eqref{cov3} satisfies $l\neq m$ and is equal to
\begin{align*}
& \quad \sum_{i=1}^n \E_H [ |U_{l,i}|^2  |U_{m,i}^2| ] f(\lambda_i)^2 + \sum_{i\neq j} \E_H [ |U_{l,i}|^2  |U_{m,j}|^2 ]  f(\lambda_i)f(\lambda_j) \\
& = \frac{\beta'}{n(\beta'n+1)} \sum_{i=1}^n  f(\lambda_i)^2 + \frac{(n-1)\beta'+1}{n(n-1)(\beta'n+1)} \sum_{i\neq j}  f(\lambda_i)f(\lambda_j)  \\
& =  \frac{\beta'}{n(\beta'n+1)}  X^\sn_1(f^2) + \frac{(n-1)\beta'+1}{n(n-1)(\beta'n+1)} (X^\sn_1(f)^2- X^\sn_1(f^2) ) \\
& = \frac{(n-1)\beta'+1}{n(n-1)(\beta'n+1)} X^\sn_1(f)^2 - \frac{1}{n(n-1)(\beta'n+1)} X^\sn_1(f^2) .
\end{align*}
And for the inner sum in \eqref{cov2} we get
\begin{align*}
& \quad \sum_{i=1}^n \E_H [ |U_{l,i}|^4  ]  f(\lambda_i)^2 + \sum_{i\neq j} \E_H [ |U_{l,i}|^2  |U_{l,j}|^2 ]  f(\lambda_i)f(\lambda_j) \\
& = \frac{1+\beta'}{n(\beta'n+1)} \sum_{i=1}^n  f(\lambda_i)^2 + \frac{\beta'}{n(\beta'n+1)} \sum_{i\neq j}  f(\lambda_i)f(\lambda_j)  \\
& =  \frac{\beta'}{n(\beta'n+1)} X^\sn_1(f)^2 + \frac{1}{n(\beta'n+1)} X^\sn_1(f^2) .
\end{align*}
So summing over $l,m$ becomes not very difficult and we obtain
\begin{align*}
& \quad \E_H \left[ X^\sn_s(f) X^\sn_t(f)\right]
\\ &= (s_n(s_n-1)+s_n(t_n-s_n)) \left( \frac{(n-1)\beta'+1}{n(n-1)(\beta'n+1)} X^\sn_1(f)^2 - \frac{1}{n(n-1)(\beta'n+1)} X^\sn_1(f^2)\right) \\
& \quad + s_n \left(\frac{\beta'}{n(\beta'n+1)} X^\sn_1(f)^2 + \frac{1}{n(\beta'n+1)} X^\sn_1(f^2) \right) .
\end{align*}
From this we have to substract the product of expectation, which we expand as
\begin{align*}
  \E_H \big[ X^\sn_s(f)\big] \E_H \big[ X^\sn_t(f) \big]
 = \frac{s_nt_n}{n^2}  X^\sn_1(f)^2 = \frac{s_n(t_n-1)}{n^2} X^\sn_1(f)^2 + \frac{s_n}{n^2} X^\sn_1(f)^2  .
\end{align*}
For the covariance we can then combine conveniently:
\begin{align} \label{varconvergence2}
& \quad \mathrm{Cov}_H(X^\sn_s(f) , X_t^\sn(f)) \notag \\
& = s_n(t_n-1) \left( \left( \frac{(n-1)\beta'+1}{n(n-1)(\beta'n+1)} - \frac{1}{n^2}\right) X^\sn_1(f)^2 - \frac{1}{n(n-1)(\beta'n+1)} X^\sn_1(f^2)\right) \notag \\
& \quad + s_n \left(\left( \frac{\beta'}{n(\beta'n+1)} - \frac{1}{n^2}\right) X^\sn_1(f)^2 + \frac{1}{n(\beta'n+1)} X^\sn_1(f^2) \right) \notag \\
& = \frac{s_n(t_n-1)}{n(n-1)} \left( \frac{1}{n(\beta'n+1)}  X^\sn_1(f)^2 - \frac{1}{\beta'n+1} X^\sn_1(f^2)\right) \notag \\
& \quad + \frac{s_n}{n} \left( \frac{-1}{n(\beta'n+1)}  X^\sn_1(f)^2 + \frac{1}{\beta'n+1} X^\sn_1(f^2) \right)  . 
\end{align}
Now $\mathbb{P}_\Lambda$-almost surely $\frac{1}{n}  X^\sn_1(f) \to \nu(f) $ and $\frac{1}{n}  X^\sn_1(f^2)\to \nu(f^2)$ by (A2), such that as $n\to \infty$, \eqref{varconvergence2} converges to
\begin{align*}
st (\beta')^{-1} ( \nu(f)^2 - \nu(f^2)) + s (\beta')^{-1} ( \nu(f^2) - \nu(f)^2)  = (s\wedge t -st ) (\beta')^{-1} ( \nu(f^2) - \nu(f)^2) ,
\end{align*}
which is precisely the right hand side of \eqref{varconvergence}. 

\bibliographystyle{alpha}
\bibliography{spectralfclt}

\begin{thebibliography}{BDMR14}

\bibitem[AGZ10]{anderson2010introduction}
G.~Anderson, A.~Guionnet, and O.~Zeitouni.
\newblock {\em An introduction to random matrices, volume 118 of Cambridge
  Studies in Advanced Mathematics}.
\newblock Cambridge University Press, Cambridge New York, 2010.

\bibitem[APS01]{albeverio20011}
S.~Albeverio, L.~Pastur, and M.~Shcherbina.
\newblock On the $1/n$ expansion for some unitary invariant ensembles of random
  matrices.
\newblock {\em Communications in Mathematical Physics}, 224(1):271--305, 2001.

\bibitem[BADG01]{arous2001aging}
G.~Ben~Arous, A.~Dembo, and A.~Guionnet.
\newblock Aging of spherical spin glasses.
\newblock {\em Probability theory and related fields}, 120(1):1--67, 2001.

\bibitem[BAG97]{arous1997large}
G.~Ben~Arous and A.~Guionnet.
\newblock Large deviations for {W}igner's law and {V}oiculescu's
  non-commutative entropy.
\newblock {\em Probability theory and related fields}, 108(4):517--542, 1997.

\bibitem[BD17]{breuer2017central}
J.~Breuer and M.~Duits.
\newblock Central limit theorems for biorthogonal ensembles and asymptotics of
  recurrence coefficients.
\newblock {\em Journal of the American Mathematical Society}, 30(1):27--66,
  2017.

\bibitem[BDMR14]{beffara2014bridges}
V.~Beffara, C.~Donati-Martin, and A.~Rouault.
\newblock Bridges and random truncations of random matrices.
\newblock {\em Random Matrices: Theory and Applications}, 3(02):1450006, 2014.

\bibitem[Bil99]{billingsley1999convergence}
P.~Billingsley.
\newblock {\em Convergence of probability measures}.
\newblock John Wiley \& Sons, 1999.

\bibitem[BLP19]{bai2019central}
Z.~Bai, H.~Li, and G.~Pan.
\newblock Central limit theorem for linear spectral statistics of large
  dimensional separable sample covariance matrices.
\newblock {\em Bernoulli}, 25(3):1838--1869, 2019.

\bibitem[BLS18]{bekerman2018clt}
F.~Bekerman, T.~Lebl{\'e}, and S.~Serfaty.
\newblock Clt for fluctuations of $\beta$-ensembles with general potential.
\newblock {\em Electronic Journal of Probability}, 23, 2018.

\bibitem[BS08]{bai2008clt}
Z.~Bai and J.~W. Silverstein.
\newblock Clt for linear spectral statistics of large-dimensional sample
  covariance matrices.
\newblock In {\em Advances In Statistics}, pages 281--333. World Scientific,
  2008.

\bibitem[BW71]{bickel1971convergence}
P.~J. Bickel and M.~J. Wichura.
\newblock Convergence criteria for multiparameter stochastic processes and some
  applications.
\newblock {\em The Annals of Mathematical Statistics}, pages 1656--1670, 1971.

\bibitem[CL95]{costin1995gaussian}
O.~Costin and J.~Lebowitz.
\newblock Gaussian fluctuation in random matrices.
\newblock {\em Physical Review Letters}, 75(1):69, 1995.

\bibitem[CS06]{collins2006integration}
B.~Collins and P.~Sniady.
\newblock Integration with respect to the {H}aar measure on unitary, orthogonal
  and symplectic group.
\newblock {\em Communications in Mathematical Physics}, 264(3):773--795, 2006.

\bibitem[{D'A}00]{d2000invariance}
A.~{D'A}ristotile.
\newblock An invariance principle for triangular arrays.
\newblock {\em Journal of Theoretical Probability}, 13(2):327--341, 2000.

\bibitem[DDN03]{d2003brownian}
A.~{D'A}ristotile, P.~Diaconis, and C.~Newman.
\newblock Brownian motion and the classical groups.
\newblock {\em Lecture Notes-Monograph Series}, pages 97--116, 2003.

\bibitem[DMR12]{donati2012truncations}
C.~Donati-Martin and A.~Rouault.
\newblock Truncations of {H}aar distributed matrices, traces and bivariate
  {B}rownian bridges.
\newblock {\em Random Matrices: Theory and Applications}, 1(01):1150007, 2012.

\bibitem[DP12]{dumitriu2012global}
I.~Dumitriu and E.~Paquette.
\newblock Global fluctuations for linear statistics of $\beta$-{J}acobi
  ensembles.
\newblock {\em Random Matrices: Theory and Applications}, 1(04):1250013, 2012.

\bibitem[Dui15]{duits2015global}
M.~Duits.
\newblock On global fluctuations for non-colliding processes.
\newblock {\em arXiv preprint arXiv:1510.08248}, 2015.

\bibitem[Dys62]{dyson1962threefold}
F.~Dyson.
\newblock The threefold way. algebraic structure of symmetry groups and
  ensembles in quantum mechanics.
\newblock {\em Journal of Mathematical Physics}, 3(6):1199--1215, 1962.

\bibitem[GNR16]{gamboa2016sum}
F.~Gamboa, J.~Nagel, and A.~Rouault.
\newblock Sum rules via large deviations.
\newblock {\em Journal of Functional Analysis}, 270(2):509--559, 2016.

\bibitem[GR11]{gamboa2011large}
F.~Gamboa and A.~Rouault.
\newblock Large deviations for random spectral measures and sum rules.
\newblock {\em Applied Mathematics Research eXpress}, 2011(2):281--307, 2011.

\bibitem[HJ77]{hoffmann1977probability}
J.~Hoffmann-J{\o}rgensen.
\newblock Probability in {B}anach space.
\newblock In {\em Ecole d’{\'e}t{\'e} de probabilit{\'e}s de Saint-Flour
  VI-1976}, pages 1--186. Springer, 1977.

\bibitem[HP06]{hiai2006semicircle}
F.~Hiai and D.~Petz.
\newblock {\em The semicircle law, free random variables and entropy}.
\newblock Number~77. American Mathematical Soc., 2006.

\bibitem[Joh97]{johansson1997random}
K.~Johansson.
\newblock On random matrices from the compact classical groups.
\newblock {\em Annals of mathematics}, pages 519--545, 1997.

\bibitem[Joh98]{johansson1998fluctuations}
K.~Johansson.
\newblock On fluctuations of eigenvalues of random {H}ermitian matrices.
\newblock {\em Duke Mathematical Journal}, 91(1):151--204, 1998.

\bibitem[Jon82]{jonsson1982some}
D.~Jonsson.
\newblock Some limit theorems for the eigenvalues of a sample covariance
  matrix.
\newblock {\em Journal of Multivariate Analysis}, 12(1):1--38, 1982.

\bibitem[KS10]{kriecherbauer2010fluctuations}
T.~Kriecherbauer and M.~Shcherbina.
\newblock Fluctuations of eigenvalues of matrix models and their applications.
\newblock {\em arXiv preprint arXiv:1003.6121}, 2010.

\bibitem[LP09a]{lytova2009central}
A.~Lytova and L.~Pastur.
\newblock Central limit theorem for linear eigenvalue statistics of random
  matrices with independent entries.
\newblock {\em The Annals of Probability}, 37(5):1778--1840, 2009.

\bibitem[LP09b]{lytova2009fluctuations}
A.~Lytova and L.~Pastur.
\newblock Fluctuations of matrix elements of regular functions of {G}aussian
  random matrices.
\newblock {\em Journal of Statistical Physics}, 134(1):147--159, 2009.

\bibitem[LP11]{lytova2011non}
A.~Lytova and L.~Pastur.
\newblock Non-{G}aussian limiting laws for the entries of regular functions of
  the {W}igner matrices.
\newblock {\em arXiv preprint arXiv:1103.2345}, 2011.

\bibitem[Meh04]{mehta2004random}
M.~L. Mehta.
\newblock {\em Random matrices}, volume 142.
\newblock Academic press, 2004.

\bibitem[M{\'S}S07]{mingo2007second}
J.~Mingo, P.~{\'S}niady, and R.~Speicher.
\newblock Second order freeness and fluctuations of random matrices: {II}.
  unitary random matrices.
\newblock {\em Advances in Mathematics}, 209(1):212--240, 2007.

\bibitem[NP10]{nourdin2010universal}
I.~Nourdin and G.~Peccati.
\newblock Universal {G}aussian fluctuations of non-{H}ermitian matrix
  ensembles: From weak convergence to almost sure {CLT}s.
\newblock {\em Alea}, 7:341--375, 2010.

\bibitem[NY16]{najim2016gaussian}
J.~Najim and J.~Yao.
\newblock Gaussian fluctuations for linear spectral statistics of large random
  covariance matrices.
\newblock {\em The Annals of Applied Probability}, 26(3):1837--1887, 2016.

\bibitem[ORS13]{o2013fluctuations}
S.~{O'R}ourke, D.~Renfrew, and A.~Soshnikov.
\newblock On fluctuations of matrix entries of regular functions of {W}igner
  matrices with non-identically distributed entries.
\newblock {\em Journal of Theoretical Probability}, 26(3):750--780, 2013.

\bibitem[Pas06]{pastur2006limiting}
L.~Pastur.
\newblock Limiting laws of linear eigenvalue statistics for {H}ermitian matrix
  models.
\newblock {\em Journal of mathematical physics}, 47(10):103303, 2006.

\bibitem[PRS12]{pizzo2012fluctuations}
A.~Pizzo, D.~Renfrew, and A.~Soshnikov.
\newblock Fluctuations of matrix entries of regular functions of {W}igner
  matrices.
\newblock {\em Journal of Statistical Physics}, 146(3):550--591, 2012.

\bibitem[Pru98]{pruss1998maximal}
A.~Pruss.
\newblock A maximal inequality for partial sums of finite exchangeable
  sequences of random variables.
\newblock {\em Proceedings of the American Mathematical Society},
  126(6):1811--1819, 1998.

\bibitem[PS11]{pastur2011eigenvalue}
L.~Pastur and M.~Shcherbina.
\newblock {\em Eigenvalue distribution of large random matrices}.
\newblock Number 171. American Mathematical Soc., 2011.

\bibitem[Rai98]{rains1998normal}
E.~Rains.
\newblock Normal limit theorems for symmetric random matrices.
\newblock {\em Probability theory and related fields}, 112(3):411--423, 1998.

\bibitem[RS06]{rider2006gaussian}
B.~Rider and J.~W. Silverstein.
\newblock Gaussian fluctuations for non-{H}ermitian random matrix ensembles.
\newblock {\em The Annals of Probability}, pages 2118--2143, 2006.

\bibitem[Shc08]{shcherbina2008invariant}
M.~Shcherbina.
\newblock Central limit theorem for linear eigenvalue statistics of
  orthogonally invariant matrix models.
\newblock {\em Journal of Mathematical Physics, Analysis, Geometry},
  4(1):171–--195, 2008.

\bibitem[Shc11]{shcherbina2011wigner}
M.~Shcherbina.
\newblock Central limit theorem for linear eigenvalue statistics of the wigner
  and sample covariance random matrices.
\newblock {\em Journal of Mathematical Physics, Analysis, Geometry},
  7(2):176–--192, 2011.

\bibitem[Sos00]{soshnikov2000central}
A.~Soshnikov.
\newblock The central limit theorem for local linear statistics in classical
  compact groups and related combinatorial identities.
\newblock {\em Annals of probability}, pages 1353--1370, 2000.

\bibitem[Sos02]{soshnikov2002gaussian}
A.~Soshnikov.
\newblock Gaussian limit for determinantal random point fields.
\newblock {\em Annals of Probability}, pages 171--187, 2002.

\bibitem[Str72]{straf1972weak}
M.~L. Straf.
\newblock Weak convergence of stochastic processes with several parameters.
\newblock In {\em Proceedings of the Sixth Berkeley Symposium on Mathematical
  Statistics and Probability}, volume~2, pages 187--221, 1972.

\end{thebibliography}

\bigskip

{\footnotesize
\noindent
TU Dortmund, \\
Fakult\"at f\"ur Mathematik, \\
Vogelpothsweg 87, 
44227 Dortmund, 
Germany, \\
jan.nagel@tu-dortmund.de
}

\end{document}